\documentclass[12pt, reqno]{amsart}
\usepackage{amssymb,amsthm,amsfonts,amsmath, amscd}

\usepackage[hypertex]{hyperref}
\usepackage{mathrsfs}

\newcommand\Y{\mathbb Y}
\newcommand\Z{\mathbb Z}

\newcommand\R{\mathbb R}

\newcommand\GT{{\mathbb{GT}}}

\newcommand\N{\mathbb N}

\newcommand\De{\Delta}

\newcommand\La{\Lambda}
\newcommand\la{\lambda}

\newcommand\om{\omega}
\newcommand\Om{\Omega}

\newcommand\wt{\widetilde}
\newcommand\wh{\widehat}

\newcommand\Tab{\operatorname{Tab}}
\newcommand\sh{\operatorname{sh}}

\newcommand\down{{\downarrow}}
\newcommand\up{{\uparrow}}
\newcommand\pd{\partial}

\newcommand\z{{(z,z')}}
\newcommand\Qz{Q^{(z,z')}_r}

\newcommand\RR{\R_{>0}}

\newtheorem{theorem}{Theorem}[section]

\newtheorem{claim}[theorem]{Claim}

\theoremstyle{definition}
\newtheorem{definition}[theorem]{Definition}
\newtheorem{remark}[theorem]{Remark}

\newtheorem{condition}[theorem]{Condition}
\newtheorem{conjecture}[theorem]{Conjecture}

\numberwithin{equation}{section}

\begin{document}

\title[Interacting particle process]{An interacting particle process\\
 related to Young tableaux}

\author{Alexei Borodin}
\address{Alexei Borodin:
\newline\indent Department of Mathematics, MIT, Cambridge, MA, USA;
\newline\indent Institute for Information Transmission Problems, Moscow, Russia}
\email{borodin@math.mit.edu}

\author{Grigori Olshanski}
\address{Grigori Olshanski:
\newline\indent Institute for Information Transmission Problems, Moscow, Russia;
\newline\indent Independent University of Moscow, Russia;
\newline\indent National Research University Higher School of Economics, Moscow, Russia}

\email{olsh2007@gmail.com}

\date{}

\begin{abstract}

We discuss a stochastic particle system consisting of a two-dimensional array
of particles living in one space dimension. The stochastic evolution bears a
certain similarity to Hammersley's process, and the particle interaction  is
governed by combinatorics of the Young tableaux.

\end{abstract}

\maketitle

\tableofcontents

\section{Introduction}

The present note is related to our papers \cite{BO-JFA12}, \cite{BO-MMJ13},
\cite{BO-2013}, and may be viewed as a companion to the latter paper.

The purpose of \cite{BO-JFA12} and \cite{BO-2013} was to construct new models
of infinite-dimensional Markov dynamics of representation-theoretic origin.

In \cite{BO-JFA12}, we studied a four-parameter family of continuous time
Feller Markov processes whose state space is the dual object $\wh{U(\infty)}$
to the infinite-dimensional unitary group $U(\infty)$. The points of the space
$\wh{U(\infty)}$ parameterize the extremal characters of $U(\infty)$.
Equivalently, the same space can also be viewed as the boundary $\pd\GT$ of the
Gelfand--Tsetlin graph.

In \cite{BO-2013}, we described a similar construction related to the infinite
symmetric group $S(\infty)$ that leads to a two-parameter family of Feller
Markov processes on an infinite-dimensional cone
$\wt\Om\subset\R^\infty\times\R^\infty\times\R$, called the {\it Thoma cone\/}.
The latter space is an amplified version of the dual object to $S(\infty)$: the
base of the Thoma cone is the {\it Thoma simplex\/} $\Om=\wh{S(\infty)}$.

The paper \cite{BO-MMJ13} explains a striking resemblance between the character
theory of the two groups, $U(\infty)$ and $S(\infty)$: we show how the
characters of the latter group arise as a degeneration of those of the former
group. In particular, we explain in what sense the Thoma cone can be viewed as
a degeneration of $\pd\GT$. Moreover, the Markov dynamics on the Thoma cone can
also be viewed as a degeneration of the Markov dynamics on $\pd\GT$.

The Markov processes from \cite{BO-JFA12} and \cite{BO-2013} may be viewed as
interacting particle processes with {\it nonlocal\/} (or {\it long-range\/})
interaction. On the other hand, as shown in \cite{BO-JFA12}, each of the Markov
processes on $\pd\GT$ is ``covered'' by a certain Markov process with {\it
local\/} interaction, living on the path space of the graph $\GT$. In the
present informal note we speculate that  similar  covering processes should
exist in the case of the Markov dynamics on the Thoma cone, too.

Although we do not have a complete proof yet, we could not resist the
temptation to introduce this model, because it provides an alternative approach
to the results of \cite{BO-2013}.

The note contains no proofs and presents a collection of conjectures and claims
based on some heuristic arguments.

\section{Generalized standard Young tableaux}

Denote $\N:=\{1,2,\dots\}$ and $\N^2:=\N\times\N$. We endow $\N^2$ with the
structure of poset such that one element  $(i,j)$ of $\N^2$ is declared to be
greater than another element $(i',j')$ if $i\ge i'$, $j\ge j'$, and at least
one of these inequalities is strict. By a (finite or infinite) Young diagram we
mean a (finite or infinite)  subset of $\N^2$ such that if it contains a given
element of the poset $\N^2$ then it also contains all smaller elements.

Let $\Y$ denote the set of all finite Young diagrams including the empty one.
We denote a Young diagram by the letter $\la$ and write $|\la|$ for the size of
$\la$. (Infinite Young diagram will appear later on.)

Recall that a {\it standard tableau\/} $\tau$ of a given shape $\la\in\Y$ is a
bijection $\tau:\la\to\{1,2,\dots,|\la|\}$ compatible with the partial order on
$\la$ inherited from $\N^2$. That is, $\tau(i,j)>\tau(i',j')$ whenever $(i,j)$
is greater than $(i',j')$.

Let $\RR$ stand for the set of strictly positive real numbers.

\begin{definition}\label{def2.A}
Following \cite{BO-MMJ13}, we define a {\it generalized\/} standard Young
tableau of a given shape $\la\in\Y$ as a function $h:\la\to\RR$ such
that
\begin{itemize}
\item $h$ is compatible with the partial order on $\la$, that is,
$h(i,j)>h(i',j')$ whenever $(i,j)$ is greater than $(i',j')$,

\item the numbers $h(i,j)$ are pairwise distinct.
\end{itemize}
\end{definition}

The conventional standard tableaux of a shape $\la\in\Y$ encode all possible
ways of building $\la$ step by step, starting with the empty set and adding an
element $(i,j)$ at each step. Likewise, a generalized standard diagram encodes
a similar process, only we imagine that elements are added at some moments
$h_1<h_2<\dots<h_{|\la|}$ on the time axis $\RR$.

\begin{definition}\label{def2.B}
Given $r\in\RR$ and $\la\in\Y$, we denote by $H_r(\la)$ the set of all
generalized standard tableaux $h$ of the shape $\la$, such that $h(i,j)<r$ for
all $(i,j)\in\la$. We agree that for the empty diagram $\varnothing$, the set
$H_r(\varnothing)$ is a singleton. Next, we set
$$
H_r:=\bigsqcup_{\la\in\Y}H_r(\la).
$$
\end{definition}

Given $h\in H_r$, we write $\sh(h)$ for the corresponding Young diagram $\la$.
It is convenient to define the ``height'' $h(i,j)$ outside $\la$, by setting
$h(i,j)=r$ for all $(i,j)\in\N^2\setminus\la$. Then we may interpret the
elements $h\in H_r$ as the ``height functions'' $h:\N^2\to(0,r]$ subject to the
following conditions:
\begin{itemize}
\item $h(i,j)<h(i,j+1)$ unless $h(i,j)=h(i,j+1)=r$,\\
\item likewise, $h(i,j)<h(i+1,j)$ unless $h(i,j)=h(i+1,j)=r$,\\
\item the set $\sh(h):=\{(i,j)\in\N^2: h(i,j)<r\}$ is finite,\\
\item the numbers $h(i,j)$, where $(i,j)$ ranges over $\sh(h)$, are pairwise
distinct.
\end{itemize}

\section{Markov dynamics on generalized tableaux}\label{sect3}

As above, we fix a number $r\in\R$. We are going to define a continuous time
Markov process on $H_r$. It depends on two continuous parameters $z$ and $z'$
subject to the following condition:
\begin{condition}\label{cond3.A}
Either both parameters $z$ and $z'$ are complex numbers with nonzero imaginary
part and $z'=\bar z$, or both parameters are real and contained in an open unit
interval of the form $(m,m+1)$ for some $m\in\Z$.
\end{condition}

This is equivalent to requiring that $(z+k)(z'+k)>0$ for every $k\in\Z$. In
particular, Condition \ref{cond3.A} implies that
\begin{equation}\label{eq3.A}
q(i,j):=(z+j-i)(z'+j-i)>0, \qquad \forall (i,j)\in\N^2,
\end{equation}
which is used in the sequel.

We need one more notation. Given $h\in H_r$ and $(i,j)\in\N^2$,  we set
\begin{equation}\label{eq3.C}
h^\down(i,j)=\begin{cases}\max\left\{h(i-1,j), h(i,j-1)\right\}, & \text{if
$i>1$
and $j>1$};\\
h(1,j-1), &\text{if $i=1$, $j>1$};\\
h(i-1,1), &\text{if $i>1$, $j=1$};\\
0, & \text{if $i=j=1$}.
\end{cases}
\end{equation}

Let us fix a couple $(z,z')$ of parameters satisfying Condition \ref{cond3.A}.
We consider an infinite collection of independent Poisson processes indexed by
the elements $(i,j)\in\N^2$, where each process is defined in the horizontal
strip
\begin{equation}\label{eq3.B}
D_r:=\{(t,x)\in\R^2, \quad -\infty<t<+\infty, \quad 0<x<r\},
\end{equation}
and the $(i,j)$th process  has constant rate equal to $q(i,j)$ (see
\eqref{eq3.A}). We denote by $\pi_r(i,j)\subset D_r$ the random point
configuration from the $(i,j)$th process. We may assume that $\pi_r(i,j)$ has
at most one point on each vertical line.

\begin{definition}\label{def3.A}
We introduce a Markov process $\wt X^\z_r$ on $H_r$, whose time evolution is
composed of two components, one deterministic and the other one stochastic:
\begin{itemize}
\item The deterministic component prescribes the height $h(i,j)$ (for every
$(i,j)\in\sh(h)$) to grow with the varying  velocity $v(y):=y(y+1)$ depending
on the current value $y:=h(i,j)$, until the height attains the maximum value
$r$, where it stabilizes.

\item The stochastic component, on the contrary,  forces the height to
instantly drop to some level. More precisely, given $(i,j)\in\N^2$, the height
$h(i,j)$ drops precisely at those time moments $t$ for which the vertical line
with abscissa $t$ contains a point $(t,x)\in \pi_r(i,j)\subset D_r$ such that
the ordinate $x$ satisfies the constraints $h^\down(i,j)<x<h(i,j)$. Then the
height instantly takes the new value equal to $x$.
\end{itemize}
\end{definition}

Note that the deterministic evolution preserves the set $H_r$: this follows
from the very definition of $H_r$ and the fact that the velocity $v(x)$ is a
strictly increasing function in $x$. Likewise, the stochastic evolution also
preserves $H_r$, by virtue of the requirement that for the drop to take place
$x$ must be greater than $h^\down(i,j)$. More pedantically, we had to include
in the above definition the additional requirement that a jump $h(i,j)\to x$ is
forbidden whenever $x$ is equal to $h(i',j')$ for some $(i',j')\in\sh(h)$
distinct from $(i,j)$, but such an event has probability 0.

We may interpret $\wt X^\z_r$ as an evolution of a system of particles on
$(0,r]$ indexed by the two-dimensional array $\N^2$. It is worth noting that
the stochastic component of the process includes an interaction between the
particles, because the lower bound $h^\down(i,j)$ for the possible drop of the
$(i,j)$th particle depends on the positions of the neighboring particles. Here
the ``neighbors'' are defined as the particles with neighboring indices
$(i-1,j)$ and $(i,j-1)$. This kind of interaction is similar to that
encountered in the exclusion processes.

On the other hand, because of the presence of two components, deterministic and
stochastic, the process $\wt X^\z_r$ may be viewed as a member of the class of
{\it piecewise deterministic Markov processes\/} introduced by Davis
\cite{Davis}.

The following claim looks very plausible to us.

\begin{conjecture}\label{conj3.A}
The process $\wt X^\z_r$ on $H_r$ does not explode and so has infinite life
time almost surely.
\end{conjecture}

\section{Evolution of Gibbs measures}

Here we discuss a connection between the process $\wt X^\z_r$ just defined and
the process $X^\z_r$ discussed in our papers \cite{BO-2013} and
\cite{BO-PTRF06}.

Let us recall the definition of $X^\z_r$ (see \cite[Definition 8.9]{BO-2013}).
This is a continuous time jump Markov process with the state space $\Y$. The
evolution is given the following $Q$-matrix (the matrix of jump rates):
\begin{equation}\label{eq4.A}
\begin{aligned}
\Qz(\la,\la+\square)&=r(z+c(\square))(z'+c(\square))
\frac{\dim(\la+\square)}{(|\la|+1)\dim\la}, \quad
\square\in\la^+,\\
\Qz(\la,\la-\square)&=(r+1)\frac{|\la|\dim(\la-\square)}{\dim\la}, \quad \square\in\la^-,\\
-\Qz(\la, \la)&=(2r+1)|\la| +rzz'.
\end{aligned}
\end{equation}
The notation is the following: $\la$ is a Young diagram; $\la\pm\Box$ is
another diagram obtained from $\la$ by appending/removing a box $\Box$;
$\la^\pm$ is the set of those boxes that can be appended to (respectively,
removed from) $\la$; $Q(\la,\la\pm\Box)$ is the rate of the jump
$\la\to\la\pm\Box$; finally, $-Q(\la,\la)$ is equal to the sum of the rates of
all possible jumps $\la\to\la\pm\Box$.

As shown in \cite{BO-2013}, this $Q$-matrix is regular, meaning that the
corresponding jump Markov process does not explode.

Let us return to the sets $H_r(\la)$ introduced in Definition \ref{def2.B}. For
every $\la\in\Y$, let $\Tab(\la)$ denote the finite set consisting of all
(conventional) standard tableaux of the shape $\la$. As seen from Definition
\ref{def3.A}, there is a natural projection $H_r(\la)\to\Tab(\la)$. For every
standard tableau $\tau\in\Tab(\la)$, the fiber of this projection over $\tau$
can be viewed as an open simplex $\De(\tau)$ of dimension $N=|\la|$, formed by
the ordered $N$-tuples of reals $(x_1,\dots,x_N)$ such that
$0<x_1<\dots<x_N<r$. It follows that $H_r(\la)$ can be viewed as a bounded open
subset of $\R^N$ whose closure $\overline{H_r(\la)}$ is a convex polytope
endowed with a triangulation.

For instance, if $\la=(2,1)$, then $N=3$, $H_r(\la)$ consists of the triples
$(x_1,x_2,x_3)=(h(1,1), h(1,2), h(2,1))$ subject to the conditions
$$
0<x_1<r, \quad 0<x_2<r, \quad 0<x_3<r, \quad x_1<x_2, \quad x_1<x_3, \quad
x_2\ne x_3,
$$
and $\overline{H_r(\la)}$ is obtained by removing the last inequality and
making the remaining inequalities weak. The polytope $\overline{H_r(\la)}$ is
the union of two closed simplices, which are singled out by the inequalities
$x_2\le x_3$ and $x_2\ge x_3$, respectively.

\begin{definition}\label{def4.A}
Recall (see Definition \ref{def2.B}) that $H_r$ is the disjoint union of the
sets $H_r(\la)$, where $\la$ ranges over $\Y$. Following \cite{BO-2013} we
define a {\it Gibbs measure\/} on $H_r$ as a probability measure such  that its
restriction to each subset $H_r(\la)$ is proportional to the Lebesgue measure.
The set of all Gibbs measures is denoted by $\mathscr G_r$.
\end{definition}

Obviously, the natural projection $H_r\to\Y$ establishes a one-to-one
correspondence between $\mathscr G_r$ and the set $\mathscr M(\Y)$ of all
probability measures on $\Y$.

\begin{claim}\label{claim4.A}
The process $\wt X^\z_r$ preserves the set\/ $\mathscr G_r$ of Gibbs measures,
and the evolution of the Gibbs measures induced by the process $\wt X^\z_r$
coincides, under the bijection $\mathscr G_r\to\mathscr M(\Y)$, with the
evolution of the probability measures on $\Y$ induced by the process $X^\z_r$.
\end{claim}

Actually, we can rigorously prove only the infinitesimal version of the claim:
the application of the infinitesimal generator of $\wt X^\z_r$ to a Gibbs
measure translates to the application of the generator of $X^\z_r$ to the
corresponding measure on $\Y$.

\begin{remark}
Here is a simple yet curious formal identity used in the proof. Assume we are
given a standard tableau of the shape $\N^2$, that is, a total order on the set
$\N^2$ compatible with its partial order. Next, for every $(i,j)\in\N^2$, set
$$
(i,j)^\down:=\begin{cases}\max\left\{(i-1,j), (i,j-1)\right\}, & \text{if $i>1$
and $j>1$},\\
(1,j-1), &\text{if $i=1$, $j>1$},\\
(i-1,1), &\text{if $i>1$, $j=1$},\\
\text{undefined}, &\text{if $(i,j)=(1,1)$},
\end{cases}
$$
and
$$
(i,j)^\up:=\min\left\{(i+1,j), (i,j+1)\right\},
$$
where the maximum and minimum are taken relative to the prescribed total order
on $\N^2$. Finally, attach to every $(i,j)\in\N^2$ a formal variable $y(i,j)$.
Then the identity in question is
$$
\sum_{(i,j)\in\N^2}\big(y((i,j)^\up)+y((i,j)^\down)-2y(i,j)\big)(z+j-i)(z'+j-i)
=2\sum_{(i,j)\in\N^2}y(i,j)
$$
with the agreement that $y((1,1)^\down):=0$.
\end{remark}

As explained in \cite{BO-2013}, the process $X^\z_r$ has a stationary
distribution, the so called {\it mixed z-measure\/}, denoted by $M^\z_r$. Let
$\wt M^\z_r$ denote the corresponding Gibbs measure on $H_r$. Claim
\ref{claim4.A} implies that $\wt M^\z_r$ serves as the stationary distribution
for the process $\wt X^\z_r$.

Claim \ref{claim4.A} says that the process $\wt X^\z_r$ in some sense
``covers'' the jump process $X^\z_r$. Note that Young diagrams $\la\in\Y$ can
be represented as particle configurations, so that $X^\z_r$, like $\wt X^\z_r$,
can also be interpreted as an interacting particle process. However, a
substantial difference between $\wt X^\z_r$ and $X^\z_r$ is that the particle
interaction is local in the former process and highly non-local in the latter
one. In this sense, $\wt X^\z_r$ seems to be simpler than $X^\z_r$.

We showed in \cite{BO-2013} that for every pair $r'>r$ there exists a ``link''
$\La^{r'}_r$ (an infinite stochastic matrix of format $\Y\times\Y$), which
intertwines the processes $X^\z_{r'}$ and $X^\z_r$ for all $\z$. For the
covering processes the picture is simpler in the sense that the processes $\wt
X^\z_{r'}$ and $X^\z_r$ are ``linked'' by an ordinary map, the truncation map
\eqref{eq5.A} that we now define.

\section{Dynamics on infinite generalized tableaux}

Given $r'>r$, we define the {\it truncation map\/} $H_{r'}\to H_r$ as the
transform
\begin{equation}\label{eq5.A}
h(i,j) \to \min(h(i,j),r), \qquad \forall(i,j)\in\N^2.
\end{equation}
Obviously, for a triple $r''>r'>r$, the composed map $H_{r''}\to H_{r'}\to H_r$
is the same as $H_{r''}\to H_r$. Therefore, we may define the projection limit
space
$$
H:=\varprojlim H_r, \qquad r\to+\infty.
$$
The elements of $H$ can be viewed as the functions
$h:\N^2\to\RR\cup\{+\infty\}$ such that
\begin{itemize}
\item $h(i,j)<h(i+1,j)$ unless $h(i,j)=h(i+1,j)=+\infty$;

\item $h(i,j)<h(i,j+1)$ unless $h(i,j)=h(i,j+1)=+\infty$;

\item the finite values $h(i,j)$ are pairwise distinct and do not have
accumulation points on $\RR$.
\end{itemize}

Note that the natural projection $H\to H_r$ is still given by \eqref{eq5.A}.
Note also that, for $h\in H$, the set
$$
\sh(h):=\{(i,j)\in\N^2: h(i,j)<+\infty\}\subseteq\N^2
$$
is a (possibly infinite) Young diagram.

Every element $h\in H$ may be interpreted as a system of particles on the
extended halfline $\RR\cup\{+\infty\}$, indexed by $\N^2$ (we assume that the
$(i,j)$th particle has coordinate $x=h(i,j)$). If $\sh(h)=\N^2$, then all
particles are on $\RR$, but it may happen that $\sh(h)$ is a proper subset of
$\N^2$; then there is an infinite reservoir of particles at infinity, indexed
by the elements of $\N^2\setminus\sh(h)$. An important requirement included in
the definition is that the particle configuration on $\RR$ is always locally
finite meaning that there are finitely many particles in any bounded interval.

\begin{claim}\label{claim5.A}
The family of Markov processes $\{\wt X^\z_r: r\in\RR\}$ is consistent with the
truncation maps $H_{r'}\to H_r$ and so there exists a unique Markov process
$\wt X^\z$ on the space $H$, consistent with this family by means of the
truncation maps $H\to H_r$.
\end{claim}

In other words, $\wt X^\z$ can be viewed as the projective limit of the
processes $\wt X^\z_r$.

The Markov dynamics of such a system, given by the process $\wt X^\z$, is
described in the same way as for the truncated processes $\wt X^\z_r$, see
Section \ref{sect3} above, only the open interval $(0,r)$ should be replaced by
the halfline $\RR$, and the endpoint $r$ is shifted to $+\infty$.

Informally, in the particle system interpretation, the deterministic component
looks as the accelerating movement of particles to the right with velocity
$v(y)=y(y+1)$, while the the concurrent stochastic component forces the
particles to instantly jump to the left. The latter component is driven by a
two-dimensional array $\{\pi(i,j): (i,j)\in\N^2\}$ of independent Poisson
processes in the half-plane $D:=\{(t,x): t\in\R, \,x>0\}$, where the $(i,j)$th
process has constant rate $q(i,j)$ given by \eqref{eq3.A}. Note that the
particles can escape from $\RR$ to infinity and return back to $\RR$.

Note also that the sample trajectories of $\wt X^\z$ look more sophisticated
than those of the truncated process $\wt X^\z_r$. This is due to a major
difference between the Poisson processes in the strip $D_r$ and in the
half-plane $D$: the pushforward of a Poisson configuration under the projection
to the $t$-axis is locally finite in the former case but not in the latter
case.

Now we extend our definition of Gibbs measures to the space $H$.

\begin{definition}\label{def5.A}
A probability measure on $H$ is said to be a {\it Gibbs measure\/} if its
pushforward under the truncation map $H\to H_r$ is a Gibbs measure on $H_r$ in
the sense of Definition \ref{def4.A}, for every $r\in\RR$. The set of all Gibbs
measures on $H$ is denoted by $\mathscr G$.
\end{definition}

The above discussion shows that the process $\wt X^\infty$ preserves the Gibbs
measures.

On the other hand, as explained in \cite{BO-MMJ13}, there is a one-to-one
correspondence $\mathscr G\leftrightarrow\mathscr M(\wt\Om)$, where $\wt\Om$ is
an infinite-dimensional cone in $\R^\infty\times\R^\infty\times\R$, called the
{\it Thoma cone\/}, and $\mathscr M(\wt\Om)$ denotes the space of probability
measures on $\wt\Om$. Further, we showed in \cite{BO-2013} that the family
$\{X^\z_r: r\in\RR\}$ of Markov processes on $\Y$ determines, via the links
$\La^{r'}_r$, a Markov process $X^\z$ on the Thoma cone $\wt\Om$. Because of
the bijection $\mathscr G\leftrightarrow\mathscr M(\wt\Om)$, this leads to the
following conclusion:

\begin{claim}\label{claim5.B}
The process $\wt X^\z$ covers the process $X^\z$ in the sense that the
evolution of arbitrary probability measures on the Thoma cone induced by $X^\z$
is the same as the evolution of the Gibbs measures on $H$ induced by $\wt
X^\z$.
\end{claim}

Here is a reformulation. Given a point $\om\in\wt\Om$, let us denote by $G_\om$
the Gibbs measure on $H$ corresponding to the Dirac measure at $\om$. The
correspondence $\om\mapsto G_\om$ can be viewed as a Markov kernel
$\La^{\wt\Om}_H=\La^{\wt\Om}_H(\om, dh)$. Then the above claim means that
$\La^{\wt\Om}_H$ intertwines $X^\z$ and $\wt X^\z$.

In \cite{BO-2013}, we showed that the process $X^\z$ on the Thoma cone has a
stationary distribution $M^\z$, which is defined by the family $\{M^\z_r:
r\in\RR\}$ via the links $\La^{r'}_r$. Let $G^\z$ stand for the corresponding
Gibbs measure on $H$; in other words, $G^\z$ is the transform of $M^\z$ by
$\La^{\wt\Om}_H$.

\begin{claim}\label{claim5.C}
The Gibbs measure $G^\z$ serves as a stationary distribution for $\wt X^\z$.
\end{claim}

\section{Concluding remarks}

Besides the truncation maps $H\to H_r$ there exist other projections under
which the Markov property of the process $\wt X^\z$ is not destroyed. Namely,
one may fix an arbitrary Young diagram, finite or infinite, and focus on those
particles that are indexed by the elements $(i,j)$  of that diagram.

The simplest such example is obtained when the diagram in question consists of
the single element $(1,1)$. Then the model represents a single particle that
moves on the extended halfline $\RR\cup\{+\infty\}$. The dynamics is described
exactly as before: a deterministic movement to the right with the velocity
$v(y)=y(y+1)$ combined with instant jumps to the left directed by the Poisson
process in the half-plane $D$ with constant rate $q(1,1)=zz'$.

A more complicated model arises when we take the infinite one-row diagram
$\{(1,j): j\in\N\}$. Then we have countably many particles, enumerated by the
numbers $j\in\N$ and occupying positions $x_1<x_2<\dots$ on $\RR$ with a
possible infinite reservoir at infinity, meaning that it may happen that
$x_j=+\infty$ for all $j$ large enough.  The dynamics in this model has a
certain resemblance with the {\it Hammersley process\/} studied in
Aldous-Diaconis \cite{AD}.

Finally, note that the process $\wt X^\z$ admits a {\it Plancherel-type\/}
degeneration in the spirit of \cite[Section 10]{BO-2013}. The idea is to scale
the space variable $x$ by substituting $x(zz')^{-1}$ instead of $x$, and pass
to a limit as both $z$ and $z'$ go to $+\infty$. Then we get a simplified
model, where the velocity in the deterministic component depends on the
coordinate $y\in\RR$ via $v(y)=y$, and the backward jump rules are the same for
all particles (they are governed by independent Poisson processes in $D$ with
constant rate $q=1$).

\section{Acknowledgements}

A.~B. was partially supported by NSF-grant DMS-1056390. G.~O. was partially
supported by a grant from Simons Foundation (Simons--IUM Fellowship), the
RFBR-CNRS grant 10-01-93114, and the project SFB 701 of Bielefeld University.

\end{document}